
\magnification=1200
\hsize=6.5truein
\vsize=9.4truein
\voffset=-.2truein
\parskip=1pc
\baselineskip=14truept
\count0=1

\font\twelverm=cmr10 at 14truept
\font\ninerm=cmr10 at 10truept
\font\svtnrm=cmr17
\font\sc=cmcsc10
\font\tenmsy=cmbx10
\font\sevenmsy=cmbx7
\font\fivemsy=cmbx5
\newfam\msyfam
\def\msy{\fam\msyfam\tenmsy}
\textfont\msyfam=\tenmsy
\scriptfont\msyfam=\sevenmsy
\scriptscriptfont\msyfam=\fivemsy

\def\CP#1{\bbc\bbp^{#1}}
\def\intro{1}
\def\const{2}
\def\des{3}
\def\app{4}

\def\F{{\cal F}}
\def\R{{\cal R}}

\def\Z{{\cal Z}}
\def\V{{\cal V}}

\def\bbc{{\msy C}}
\def\bbp{{\msy P}}
\def\bbz{{\msy Z}}

\def\grF{\Phi}
\def\grQ{\Psi}

\def\grL{\Lambda}
\def\gre{\epsilon}
\def\grf{\phi}
\def\grg{\gamma}
\def\grl{\lambda}
\def\grm{\mu}
\def\grn{\nu}
\def\grq{\psi}
\def\grp{\pi}
\def\grs{\sigma}
\def\grx{\xi}

\outer\def\abstract #1\par{\centerline{\vbox{\hsize = 5.85truein
          \baselineskip = 12.5truept
          \ninerm
          {\sc Abstract}: #1 }}}
\outer\def\section #1. #2\par{\noindent{\bf \S #1. #2}\par\medskip}
\def\endsection{\bigskip}  
\outer\def\proclaim #1: #2\par{\medskip
    \noindent{\sc#1:\enspace}{\tensl#2}\par
    \ifdim\lastskip<\medskipamount \removelastskip\penalty55
    \medskip\fi}
\def\blbx{\vrule height6pt width4pt depth1pt}
\def\fract#1#2{\raise4pt\hbox{$ #1 \atop #2 $}}
\def\decdnar#1{\phantom{\hbox{$\scriptstyle{#1}$}}
\left\downarrow\vbox{\vskip15pt\hbox{$\scriptstyle{#1}$}}\right.}
\outer\def\endproof{\hfill $\blbx$\medskip}
\outer\def\proof{\noindent{\sc Proof:\enspace}}
\def\ra#1{\hbox to #1pc{\rightarrowfill}}

\bigskip
\centerline{\svtnrm The Space of Harmonic Maps}
\centerline{\svtnrm from the 2-sphere to the Complex Projective Plane}
\bigskip
\centerline{T. Arleigh Crawford}
\medskip
\centerline{Fields Institute for Research in Mathematical Sciences}
\centerline{Toronto, Canada}
\bigskip

\abstract In this paper we study the topology of the space
of harmonic maps from~$\scriptstyle S^2$ to~$\scriptstyle \CP 2$.  We
prove that the subspaces consisting of maps of a fixed degree and
energy are path connected.  By a result of Guest and Ohnita it follows
that the same is true for the space of harmonic maps to~$\scriptstyle
\CP n$ for~$\scriptstyle n\geq 2$.  We show that the components of
maps to~$\scriptstyle \CP 2 $ are complex manifolds.

\section \intro. Introduction

Harmonic maps from the Riemann sphere to complex projective space
are critical points of the energy functional
$$E:C^\infty(S^2,\CP n)~\ra{1.5}~[0,\infty)$$
defined on the space of smooth maps.  As solutions of a classical
variational problem harmonic maps have been studied extensively,
especially with regard to questions of existence, uniqueness,
regularity, etc.  It is now well known that all of the harmonic maps
from the Riemann sphere to complex projective space can be constructed
from holomorphic maps.  In this paper we study some global topological
properties of the solution space.

Holomorphic (and anti-holomorphic) maps are the absolute minima of $E$
in each path component of $C^\infty(S^2,\CP n)$.  A holomorphic
map~$f:S^2\ra{1}\CP 2$ is full if its image
is not contained in any  projective line.  The
following theorem is the starting point for this paper:

\proclaim Theorem \intro.1: [EW2] The set of all non-minimal harmonic maps
$\grf:S^2\ra{1}\CP 2$ is in 1-1 correspondence with the set of full
holomorphic maps~$f:S^2\ra{1}\CP 2$.

This is, in fact, a specific case of a more general theorem which
describes how to construct all harmonic maps to $\CP n$ set
theoretically.  It is due originally to Din and Zakrzewski [DZ].
The paper [EW2] of Eells and Wood gives an excellent description
for a mathematical audience and we refer to the construction as the
Eells-Wood construction.  The reader is
also directed to [Bu], [G], [La] for other descriptions.
With the Eells-Wood construction in hand it is possible to study the
global, topological properties of the space of harmonic maps.
Much is, in fact already known.  Holomorphic and anti-holomorphic
maps are local minima of~$E$ and the topology of these components is
very well understood (see [S], [CCMM], [CS], [H]).
In [C] we studied the subspaces of full holomorphic maps.  Some results
have also been obtained for spaces of harmonic maps of $S^2$ into
other Riemannian manifolds (see [A], [FGKO], [K], [V], [Lo]).

Let $Harm(\CP n)$ be the space of harmonic maps~$S^2\ra{1}\CP n$
with the topology it inherits as a subset of the
space~$Map(S^2,\CP n)$ of all
continuous maps.
$Map(S^2,\CP n)$ has connected components, $Map_k(S^2,\CP n)$, indexed by
the degree~$k$ of the individual maps.
The critical values of~$E$ are discrete so $Harm(\CP n)$ is a
disjoint union of the sets $Harm_{k,E}(\CP n)$ consisting of maps of
degree~$k$ and energy~$E$.  Let $Hol_k(\CP n)$ denote the space
of holomorphic maps of degree~$k$ when~$k\geq 0$ and anti-holomorphic
maps of degree~$k$ when~$k<0$.  We normalize the energy functional so
that $E(f)=k$ for a holomorphic map~$f$ of degree~$k$.

Guest and Ohnita [GO] have conjectured that the spaces
$Harm_{k,E}(\CP n)$ are connected.  They show that,
for~$n\geq 3$, any harmonic map $\grf:S^2\ra{1}\CP n$ can be continuously
deformed through harmonic maps to a map whose image lies
in~$\CP{n-1}\subset\CP n$.  Thus to prove the
conjecture it suffices to prove it in the case where~$n=2.$  The main
result of this paper is the following theorem, which affirms this
conjecture:

\proclaim Theorem \intro.2: The path components of~$Harm(\CP 2)$
are the minimal sets, $Hol_k(\CP 2)$, and the non-minimal critical
sets, $Harm_{k,E_r}(\CP 2)$, where the critical energy values
are $E_r = 3|k|+2r+4$ indexed by a non-negative integer~$r$.
Moreover these non-minimal components can
be given the structure of complex manifolds of complex
dimension~$3|k|+r+8$.

\proclaim Remark:  In a recent preprint [LW] Lemaire
and Wood show that $Harm_{k,E_r}(\CP 2)$ is, in fact,
a smooth submanifold of the space of smooth maps
for~$r\leq {1\over 2}(k+1)$.

The proof follows from an examination of the Eells-Wood construction.
The set-theoretic assignment of Theorem~\intro.1 cannot be
continuous.  To see this recall that a map~$f:S^2\ra{1}\CP 2$
is ramified at~$x\in S^2$ if $df_x=0.$  If $f$ is
holomorphic and non-constant the set of points at which $f$ is
ramified is finite and the ramification index of~$f$ is the
number of ramification points, counting multiplicities.  Let~$r$
be the ramification index of a map $f$, then~$r\leq 2k-2$ and,
unless the image of~$f$ lies in a projective line,
$r\leq{1\over 2}(3k-6)$.  The harmonic map
corresponding to~$f$, in Theorem~\intro.1, has
degree~$k-2-r$ and energy~$3k-2-r$.  However there are
plenty of maps with the same degree and different ramification
indices so the connected space of all full maps
of degree~$k$ is mapped by this correspondence to a number of
different, disjoint pieces of $Harm(\CP 2)$.

The first result of this paper is the following lemma which says that
these are the only
discontinuities.  Let $Hol_{k,r}(\CP 2)\subset Hol_k(\CP 2)$ be
the subspace of maps with ramification index~$r$.  Using the
Eells-Wood construction we prove the following:

\proclaim Lemma \intro.3: For~$0\leq r\leq k-2$ there is a
homeomorphism
$$\grF_{k,r}:Harm_{k-2-r,3k-2-r}(\CP 2)\cong Hol_{k,r}(\CP 2)$$

To get the results about connectedness and smoothness we prove:

\proclaim Theorem \intro.4: For $r\leq k-2$, $Hol_{k,r}(\CP 2)$ is
smooth connected complex submanifold of
$Hol_k(\CP 2)$ of complex dimension~$3k - 2r + 2$.

For the non-negative degree components Theorem~\intro.2 follows
from Theorem~\intro.4 and Lemma~\intro.3.
Using the fact that the involution $z\mapsto\bar z$ induces a homeomorphism
$$Harm_k(\CP 2)~\cong~Harm_{-k}(\CP 2)$$
which preserves energy we obtain the result for all cases.

The paper is organized as follows: in section two we recall the
construction of harmonic maps from holomorphic maps which lies
at the heart of the proof of Theorem~\intro.1.  We prove that,
when restricted to~$Hol_{k,r}(\CP 2)$, this construction produces a
continuous map and prove Lemma~\intro.3.  In section three we describe
the geometry of the spaces~$Hol_{k,r}(\CP 2)$ and prove Theorem~\intro.4.
In the final section we give concrete descriptions of some of the
simpler components and, where possible, compute their cohomology groups.

The research for this paper was conducted at the University of New
Mexico while I was a Ph.D. candidate and at McGill
University where I held an NSERC Post-Doctoral
Fellowship.  The final version of this paper was prepared while
I was a Post-Doctoral Fellow at the Fields Institute
for Research in Mathematical Sciences and a temporary visitor at
the University of Toronto.
I would like to thank all these institutions for their hospitality and
support.  I would like to thank my thesis advisor, Ben Mann, for
suggesting the problem and for his guidance.  I also benefitted from
conversations with Jacques Hurtubise and Martin Guest.  I am endebted
to two referees of earlier versions of this paper for extensive
and useful comments and corrections.  Finally I would like to thank
John Wood for pointing out a serious gap in the proof of an earlier
version of Theorem~\intro.4.

\endsection

\section \const. The Construction of Harmonic Maps

In this section we describe the Eells-Wood construction and show
that it restricts to a continuous, proper map on the subsets of
holomorphic maps with fixed degree and ramification index.
A number of descriptions of this construction exist in the
literature.  The one we use below is closest in spirit to [Bu] or
[EW2].  For brevity let $Hol_k = Hol_k(\CP 2)$ and
$Hol_{k,r}=Hol_{k,r}(\CP 2)$.

A holomorphic map~$f:S^2\ra{1}\CP 2$ may be defined by letting
$$f(z) ~=~ [p_0(z),p_1(z),p_2(z)],$$
where $z$ is a complex coordinate
on~$\bbc\cong S^2\setminus\{\infty\}$ and $[u_0,u_1,u_2]$ are
homogeneous coordinates on~$\CP 2$.  The $p_i$ are polynomials
which share no common zero.
The degree of~$f$ is the maximum of the degrees of the~$p_i$.
Taking coefficients of the $p_i$ as homogeneous coordinates
gives an embedding $Hol_k\subset\CP N$ as an open submanifold,
where $N=3k+2$.  Let $p:\bbc\ra{1}\bbc^3$ be the polynomial function
$$p(z) ~=~ (p_0(z),p_1(z),p_2(z)).$$
This is just a lift of~$f$ to~$\bbc^3$ over the coordinate
patch.
We will often write $[p_0,p_1(z),p_2]$ or even $[p]$ for~$f$.

If $f\in Hol_k$ then $p(z)$
and $p'(z)$ will be linearly independent for all but a
finite number of points.  The map $h=p\wedge p'$, given by
$$h(z) ~=~ p(z)\wedge p'(z)\in{\bigwedge}^2\bbc^3,$$
is also polynomial. That is, identifying
$\bigwedge^2\bbc^3\cong\bbc^3$, we can write
$$h(z) ~=~ (h_0(z),h_1(z),h_2(z)),$$
where the $h_i$ are polynomials of degree less than or equal to
$2k-2$.  If $f$ is unramified then the $h_i$ have no common zeros
and $[h]$ is holomorphic map to $\bbp\left(\bigwedge^2\bbc^3\right)
\cong G_2(\bbc^3)$ of degree~$2k-2$.  Here $G_2(\bbc^3)$ denotes the
Grassmanian of 2-planes in $\bbc^3$.  If
$f$ is ramified at $z$ then $h(z)=0$, and if $f$ is ramified at
$\infty$ then the $h_i$ will have degree strictly less than~$2k-2$.
We may write $h_i = bq_i$, for $i=1,2,3$, where $b$ is a greatest
common divisor of the~$q_i$.  Let $2k-2-r$ be the maximum of
the degrees of the~$q_i$ and let $q(z) =
(q_0(z),q_1(z),q_2(z))$, then $f_1 = [q]$ is a well-defined holomorphic
map to $G_2(\bbc^3)$ of degree~$2k-2-r$.  The integer $r$ is the
ramification index of~$f$.  The map $f_1$ is called the first
associated curve of $f$.

The line $f(z)$ is contained in the
plane $f_1(z)$ with complex codimension 1.  Thus we can define a map
$$\grf_s:S^2~\ra{1.5}~\CP 2$$
by taking
$$\grf_1(z) ~=~ f_1(z)\cap f^\perp(z).$$
The map $\grf_1$ is harmonic and the assignment
$f\mapsto \grf_1$ is the correspondence of Theorem~\intro.1.
This assignment, restricted to a fixed degree~$k$ and
ramification index $r$, defines the map
$$\grF_{k,r}:Hol_{k,r}~\ra{1.5}~
Harm_{k-2-r,3k-2-r}(\CP 2).$$
We will prove Lemma~\intro.3 by showing that $\grF_{k,r}$ is continuous
and proper.

Let $\V_d\subset\bbc[z]$ be the subspace of polynomials of degree
less than or equal to~$d$.  We can stratify the projective
space~$\bbp \V_k^3$ by taking
the subsets~$S_r$ of points $[p_0,p_1,p_2]\in\bbp \V_k^3$ such that
if~$b$ is a greatest common divisor of the~$p_i$, and we
write $p_i=bq_i$, then $k-r$ is the maximum of the degrees of the $q_i$.
Note that $S_0\cong Hol_k$.  In fact, the assignment
$$([b],[q_0,\ldots,q_n])~\mapsto~[bq_0,\ldots,bq_n]$$
defines an embedding
$$\grx:\bbp \V_r\times Hol_{k-r} \to\bbp \V_k^3$$
which shows that for $0\leq r<k$, $S_r$ is a submanifold
of~$\bbp \V_k^3$.  Note also that the closure of $S_r$ is contained
in the union of the strata $S_{r'}$ for $r'\geq r$.

Note that $Hol_{k,r}$ is just the
inverse image of $S_r$ under the map
$$\grQ:Hol_k~\ra{1.5}~\bbp \V^3_{2k-2}$$
given by $[p]\mapsto [p\wedge p^\prime]$ where $p$ is a
3-tuple of polynomials.  It follows that the first
associated curve $f_1\in Hol_{2k-2-r}(G_2(\bbc^3))$
depends continuously on~$f=[p]$.  The remainder of the
construction is manifestly continuous and it follows that
$\grF_{k,r}$ is continuous.
We also remark that the first factor,~$[b]\in\bbp \V_r$,
of~$\grx^{-1}(\grQ(f))$ also depends continuously on~$f$.
This is the ramification divisor of~$f$ and we denote it by~$\R(f)$.

\proclaim Lemma \const.2:  $\grF_{k,r}$ is proper.

\proof  The proof follows the proof of Lemma 3.3 in [FGKO].
Suppose we have a sequence $\{\grf_n\}$ converging to $\grf$ in
$\grF_{k,r}(Hol_{k,r})\subset Harm_{k-2-r,3k-2-r}(\CP 2)$
and suppose $\{f_n\}\subset Hol_{k,r}$ is such
that $\grf_n = \grF_{k,r}(f_n)$ for each $n\geq 0$.  It will suffice
to find a convergent subsequence of $\{f_n\}$.  We have $Hol_{k,r}
\subset\bbp \V_k^3$.  Since the latter space is
compact there is a subsequence which converges to a point in
$[p]\in\bbp \V_k^3$.  We must show that $[p]\in Hol_{k,r}$.
Suppose $[p]$ is in the stratum $S_m$ for
$m\geq 0$.  Then we can write $[p] = [bq_0,bq_1,bq_2]$ where the
$q_i$ have no common zero. Thus $[q]=[q_0,q_1,q_2]\in Hol_{(k-m)}$.

Similarly, by choosing a further subsequence if necessary, we can
assume that $\grQ(f_n)$ converges to some point
$[s]\in\bbp\V_{2k-2}^3$.  Since $\grQ(f_n)\in S_r$ we must have
$[s]\in S_{r'}$ with $r'\geq r$.  Write $[s] = [dt_0,dt_1,dt_2]$
with the $t_i$ coprime. So $[t]=[t_0,t_1,t_2]$ is in
$Hol_{2k-2-r'}(G_2(\bbc^{n+1}))$.

Let ${\cal Z}\subset S^2$ be the set of zeros of $b$ and $d$ and include
the point at infinity if either $\deg b<m$ or $\deg d<r'$.  The line
$[q(z)]$ is contained in the plane represented by $[t(z)]$ for all
$z$ so let
$$\grq(z)~=~ [t(z)]\cap [q(z)]^\perp.$$
$\grq$ is a harmonic map which agrees
with $\grf$ on $S^2\setminus {\cal Z}$.  So, by the unique
continuation property of harmonic maps, we must
have~$\grq=\grf$.  But
$$\eqalign{E(\grq) &~=~(2k-2-r') + (k-m)~=~3k-2-r'-m\cr
           \deg \grq &~=~(2k-2-r') - (k-m)~=~k-2-r'+m.\cr}$$
Requiring $E(\grq)=E(\grf)$ and $\deg\grq = \deg\grf$ we must have
$m=0$ and $r'=r$.  Thus $[p]\in Hol_{k,r}$.\endproof

\endsection

\section \des. The Desingularizing Variety

In this section we study the geometry of the strata
$Hol_{k,r}$.  We start by construction a filtration
$$Hol_k~=~\F_0\supset\F_1\supset\cdots\supset\emptyset$$
by closed subsets.  The strata are the differences of successive
elements in this filtration.  We construct varieties which sit over
the $\F_r$ and show that these varieties are smooth.
This is sufficient to show that $Hol_{k,r}$ is smooth and
connected.

To define the filtration let
$$\F_r ~=~ \bigl\{f\in Hol_k\bigm| \hbox{$f$ has ramification
                                index $\ge r$}\bigr\}.$$
Then ~$Hol_{k,r}=\F_r\setminus\F_{r+1}$.
It is useful to think of the ramification
divisor~$\R(f)\in \bbp \V_r$ as being in the symmetric product~$SP^r(S^2)
=(S^2)^r/{\cal S}_r$ where ${\cal S}_r$ is the symmetric group on
$r$~letters.  An explicit homeomorphism~$SP^r(S^2)\cong \bbp \V_r$
is given by mapping an unordered $r$-tuple~$\langle x_1,\ldots,
x_r\rangle$ to the equivalence class of a polynomial whose zeros
are precisely those points~$x_i$ which are not equal to~$\infty$.
We will say that
$\langle x_1,\ldots,x_s\rangle\in SP^s(S^2)$
divides~$\R(f)$ if~$\R(f) =
\langle x_1,\ldots,x_s,y_{s+1},\ldots,y_r\rangle$ for
some~$y_{r+1},\ldots,y_r$.  If all the points~$x_i$ are finite
then this is just the usual notion of polynomial division.

Let
$$X_r~=~\{([a],f)\in\bbp \V_r\times Hol_k\mid
               \hbox{$[a]$ divides $\R(f)$}\}.$$
By projecting onto the second factor we get a quotient
map $p_r: X_r\ra{1}\F_r$.  The inverse image~$p_r^{-1}(f)$
counts the (finite) number of elements $[a]$ wich
divide~$\R(f)$.  For maps with ramification index
exactly $r$ there is only one point in the inverse image and~$p_r$
restricts to a homeomorphism
$X_r\setminus p_r^{-1}(\F_{r+1})\cong
\F_r\setminus\F_{r+1} = Hol_{k,r}$.
We will prove the following:

\proclaim Lemma \des.1: For $r\leq k-2$ the
spaces $X_r$ are  path-connected complex submanifolds of
$Hol_k$ of complex codimension~$2r$.

This will imply Theorem \intro.4.  First of all, it identifies
$Hol_{k,r}$ with an open submanifold of a complex manifold
of the correct dimension.
Second, $Hol_{k,r}$ is connected since $p_r^{-1}(\F_{r+1})$ is
a proper algebraic subset and cannot disconnect a smooth variety.

In order to study the geometry of $X_r$ we need to characterize
the condition that $[a]$ divide~$\R(f)$.  Let~$f=[p_0,p_1,p_2]$.
Recall that in~\S\const\ we saw
that we could write $\grQ(f) = [bq_0,bq_1,bq_2]$ where $[b]=\R(f).$
The polynomial factors of $\grQ(f)$ are of the form
$p_ip_j' - p_i'p_j,$
for~$i < j$.  If $\deg a = r$ then $[a]$ divides
$\R(f)$ if $a$ divides $p_ip_j' - p_i'p_j$ for all
$0\leq i<j\leq 2$. These conditions are not independent:  Suppose
$p_0$~and~$a$ are coprime and that $a$~divides
$p_0p_i' - p_0'p_i$, for $i=1,2$.  Now
$$p_1(p_0p_2'-p_0'p_2)-p_2(p_0p_1'-p_0'p_1) = p_0(p_1p_2'-p_1'p_2)$$
and if $a$ divides both terms on the lefthand side it must
also divide $p_1p_2'-p_1'p_2$.

Let $X'_r\subset X_r$ be the subset of pairs~$([a],[p_0,p_1,p_2])$
such that $\deg a = r$, and $a$~and~$p_0$ are coprime.  Lemma~\des.1
will follow from the next two lemmas.

\proclaim Lemma \des.2: Every point in $X_r$ is contained in a
neighbourhood homeomorphic to $X'_r$.

\proof By a change of complex coordinate on $S^2$ we may assume
that the configuration associated to~$[a]$ does not
include the point at infinity, so~$\deg a = r$.

Now $\bbp GL(3,\bbc)$ acts on $\CP 2$ by complex,
linear bi-holomorphisms.  Thus it acts by composition
on $Hol_k$ leaving the subspaces $Hol_{k,r}$
invariant.  In fact, for $A\in\bbp GL(3,\bbc)$, $\R(A\cdot f) =
\R(f)$.
Write $f=[p_0,p_1,p_2]$.  It suffices to find $\gre_0,\gre_1,
\gre_2\in\bbc$ so that $\gre_0p_0+\gre_1p_1 +\gre_2p_2$
is prime to $a$.  Since
no zero of $a$ can be a zero of all the $p_i$ this condition
is satisfied by a generic choice of $\gre_i$. \endproof

\proclaim Lemma \des.3: For $r\leq k-2$ the spaces
$X'_r$ are complex manifolds.

\proof We will prove the lemma by giving an explicit description
of the space as a smooth pull-back.
Let $\V^+_d$ denote the set of monic polynomials in $\V_d$.  Let $\Z$
be the set of pairs $(a,p)\in \V^+_r\times \V_k$ such that $a$ and $p$
are coprime.  Let $Mat(s,t)$
be the space of $s\times t$ complex matrices. For~$s\leq t$,
let $Mat^*(s,t)$ be the open subset of matrices with rank~$s$.
We may identify
$Mat^*(s,t)$ with a fattened version of the Stiefel manifold
$V_s(\bbc^t)$ of $s$-frames in $\bbc^t$ by thinking of the
$s$~linearly independent rows of a matrix as a frame.
Now define a map
$L:\Z \ra{1.5} Mat(r,k+1)$
as follows:  Given $(a,p)\in \Z$ we can construct a linear map
$$L(a,p)\in Hom_\bbc(\V_k,\V_{r-1})\cong Mat(r,k+1)$$
by $L(a,p)\cdot u = [pu'-p'u]_a$ where, for any polynomial $q$, $[q]_a$
is the congruence class of $q$ mod $a$.

Next consider the space
$$E~=~\{(A;u_1,u_2,u_3))\in Mat(r,k+1)\times V_3(\bbc^{k+1})
      \mid u_i\in\ker A,i = 1,2,3\}.$$
Projection onto the second factor $E\ra{1}V_3(\bbc^{k+1})$
makes $E$ a vector bundle.  To see this first note that the
condition that each of the $u_i$ be in $\ker A$ is equivalent to
requiring that each of the rows of $A$ be in the kernel of the
matrix with rows $u_1$, $u_2$ and $u_3$.
So, if we map
$$V_3(\bbc^{k+1})~\ra{1.5}~G_{k-2}(\bbc^{k+1})$$
by associating
to each 3-frame a $3\times (k+1)$ matrix which we map to its
kernel, then $E$ is the pullback of the $r$-fold Whitney sum
of the canonical $\bbc^{k+1}$-bundle over $G_{k-2}(\bbc^{k+1})$.

We are now in a position to describe $X'_r$.  Consider the
diagram
$$\matrix{&&E\cr
          &&\decdnar{\grq}\cr
          \Z & \fract{\grf}{\ra{1.5}}
           & Mat(r,k+1)\times \bbc^{k+1}\cr}$$
where $\grq$ is the projection $(A;u_1,u_2,u_3)\mapsto (A,u_1)$
and $\grf$ sends $(a,p)\mapsto (L(a,p),p)$.  The pullback of this diagram
can be described as the set of points $(a,p_0,p_1,p_2)\in
\V^+_r\times \V_k^3$ with $(a,p_0)\in\Z$ and $p_i\in\ker L(a,p_0)$
for $i = 1,2$.  It is clear that $d\grq$ maps onto $TMat(r,k+1)$
and $d\grf$ maps onto $T\bbc^{k+1}$ so $\grq$ and $\grf$ are transversal
and the pullback is a manifold.  Finally we projectivise by identifying
$(a,p_0,p_1,p_2)\sim (a,\grl p_0,\grl p_1,\grl p_2)$.  Then we can equate
$X'_r$ with the open submanifold comprised of equivalence classes for which
$p_0,p_1,p_2$ have no common zeros so that they define a holomorphic map.
\endproof

\endsection
\vfill
\eject

\section \app. Appendix: Some Examples

With the description of the strata $Hol_{k,r}$ developed in the
previous section it is possible to obtain explicit geometrical models
for the first few strata.  In this section we describe some
examples.  Let $G = GL(3,\bbc)/Z$ where $Z \cong \bbc^*$ is the centre.
This is the automorphism group of $\CP 3$ and so it acts on
$Hol_k$.  The condition that the image not lie in any
$\CP 1$ means that $G$ acts freely
on $Hol_k$ and it fixes
the strata $Hol_{k,r}$.  Thus the components $Harm_{D,E}(\CP 2)$
have free $G$-actions.  The construction in the last section can be
modified slightly to describe these components as pull-backs of
a canonical principle~$G$-bundle.  We use the machinery of algebraic
topology to make calculations of some of the cohomology groups.  This
section assumes some background in algebraic topology.

The first non-trivial case is degree~$2$.  Consider a triple of
degree~$2$ polynomials, $(p_0,p_1,p_2)$.  In order to define a full map
the~$p_i$ must be linearly independent in the space of polynomials,
$\V_2$.  This condition
would be violated if the $p_i$ had a common zero.  For the same
reason the map they define must be unramified.  We may use the
coefficients of the three polynomials to form a matrix in $GL(3,\bbc)$,
defined up to multiplication by a non-zero scalar.
Thus $Hol_{3,0}\cong G$.  This
result, for based maps, appears in [C].

The group
$G$ has the homotopy type of $\bbp U(3)=U(3)/S^1$.
The cohomology of $\bbp U(3)$ is known [BB].  For $p\neq 3$
$$H^*(\bbp U(3);\bbz/p) = H^*(SU(3);\bbz/p) = \grL[e_3,e_5].$$
For $p=3$
$$H^*(\bbp U(3);\bbz/3) = \grL[e_1,e_3]\otimes\bbz/3[x_2]/x_2^3=0.$$
Where $\grL[a_i,...]$ and $\bbz/3[b_i,...]$ denote an exterior algebra
and a polynomial algebra on the given generators.  The subscripts denote
the dimensions of the generators.

In degree~3 we may again describe a full map by three linearly
independent polynomials.  Consider the condition that
$$\grm p_i(z) = \grl p'_i(z),\quad i=1,2,3\leqno(\app.1)$$
for some $z\in\bbc$.  If
$\grl\neq 0$ this corresponds to ramification at $z$.  If $\grl=0$ then
the $p_i$ all vanish at $z$.  In this way we see the possible ramifications
at $z$ being parameterized by $\grm\in\bbc$ with the extreme
case,~$\grm=\infty$, corresponding to the simultaneous vanishing of all three
polynomials.

To see what happens as $z\rightarrow\infty$ we change
coordinates to $\grx = z^{-1}$.  We look at new polynomials~$q_i$ defined
by requiring $q_i(\grx) = \grx^3 p_i(\grx^{-1})$ for $\grx\neq 0$.
And a new condition for ramification corresponding to \app.1:
$\grn q_i(\grx) = \grg q'_i(\grx)$.  Gluing these two pictures
together along the overlap $\bbc\setminus\{0\}$ we obtain an bundle over
$S^2$ with fibre $S^2$.  We denote the total space of this bundle by
$X$.  It is useful to think of $X$ as a line bundle,
$Y$, compactified by adding a section at infinity.  The finite part of
the fibre over $z$ gives the data for ramifications at $z$,
and the extra point at infinity corresponds to the condition that all
three polynomials vanish at $z$.  An explicit
calculation of the transition functions shows that $c_1(Y) = 2$.

To obtain descriptions of the two strata in $Hol_3$ we use the
ramification data to construct a pullback bundle.  This is essentially
the same construction used in the last section.  Matters are considerably
simplified by the fact that, in degree~3, $L$ maps $\Z$ into $Mat^*(r,4)$
the full-rank matrices.  So the kernel of $L(a,p)$ has constant
dimension and $E$ is the pull back over
$\ker:Mat^*(r,4)\ra{1}G_{4-r}(\bbc^4)$ of a principle $GL(3,\bbc)$-bundle.
Another way of putting this is that
condition \app.1 is always non-degenerate and
defines a 3-dimensional subspace of the space of polynomials
$\V_3\cong\bbc^4$. This
extends over the fibre at infinity to give a map
$\grf:X\hookrightarrow G_3(\V_3)\cong G_3(\bbc^4)\cong\CP 3$.

Since $Y$ parameterizes the ramification data we can write
$$Hol_{3,1}=\{(y,[p_0,p_1,p_2])\in Y\times\bbp \V_3^3\mid
                       \hbox{$p_i$ are linearly independent and
                        span $\grf(y)$}\}.$$
Let $V_3(\bbc^4)$ be the Stiefel manifold of $3$~frames in~$\bbc^4$
and let~$\bbp V_3(\bbc^4) = V_3(\bbc^4)/\bbc^*$.  Then the canonical
$GL(3,\bbc)$-bundle projection is $C^*$-invariant and we obtain a
principle $G$-bundle
$$\grp:\bbp V_3(\bbc^4)\rightarrow G_3(\bbc^4)\cong \CP 3$$
And so $Hol_{3,1}$ is the total space of the
pullback of $\grp$ over
$$Y~\subset~X~\fract{\grf}{\ra{1.5}}~\CP 3.$$
By Corollary \intro.4 this gives us a description of $Harm_{0,6}(\CP 2)$.

We can compute the cohomology of this space.  A straightforward
calculation shows that $\grf$ restricted to the zero section in
$Y$ has degree~3 and so the first Chern class pulls back to
three times a generator in $H^2(Y)$.
We need to know the cohomology of $\bbp V_3(\bbc^4)$ and the differentials
in the Serre spectral sequence for the bundle $\grp$.  The cohomology
calculation is a straightforward application of Baum's results regarding
the Eilenberg-Moore spectral sequence for the cohomology of a homogeneous
space [Ba] and from this we can deduce the
neccesary differentials.  The $E_2$ term is
$$E_2(\grp)^{*,*} = H^*(\CP 3;\bbz/p)\otimes H^*(\bbp U(3);\bbz/p).$$
We can describe the differentials in terms of the generators for
the cohomology of the fibre given above.  Let $b$ be the mod~$p$
reduction of the first Chern class.  There are three cases:
For $p=2$ there are no non-trivial differentials.  For $p=3$ there
is only one non-trivial differential generated by $d_2(e_1) = b$.
For $p>3$ the only non-trivial differentials are generated by
$d_4(e_3) = b^2$.  This is sufficient to determine the differentials in
the spectral sequence for cohomology with $\bbz/p$ coefficients
of the pullback bundle
$$\bbp U(3)~\ra{1.5}~Hol_{3,1}~\ra{1.5}~Y.$$
Since the base has the homotopy type of~$S^2$ the only possible
differentials are at $E_2$.  For $p\ne 3$ there are no non-trivial
differentials in $E_2(\grp)$.  For $p=3$, $b$ pulls back to zero so
the only possible differential is zero.
In either case the spectral sequence collapses at $E_2$ and
$$H^*(Harm_{0,6}(\CP 2);\bbz/p)\cong
               H^*(S^2;\bbz/p)\otimes H^*(\bbp U(3)\bbz/p).$$
The spectral sequence does not completely determine the cup products
so this need not be an algebra isomorphism.

We can also describe $Hol_{3,0}$.  It is the restriction of
the bundle $\grp$ to $\CP 3\setminus\grf(X)$.  This gives a geometric
description of $Harm_{1,7}(\CP 2)$, the first non-minimal critical
level in the degree~1 component, as a principle $G$-bundle
over the compliment of $X$ in $\CP 3$.  This compliment is not
simply connected and so using the Serre spectral sequence to
compute cohomology groups will involve understanding the system
of local coefficients.

\endsection

\bigskip
\centerline{\twelverm References}
\bigskip
\hfuzz=4pt

\item{[A]} {\sc C.K. Anand}, {\tensl Uniton Bundles}, to appear Comm.
Anal. Geom.

\item{[BHMM]} {\sc C.P. Boyer, J.C. Hurtubise, B.M. Mann, R.J. Milgram},
{\tensl The topology of the instanton moduli spaces.  I: The
Atiyah-Jones Conjecture},  Ann. Math. (2) 197(1993), 561--609.

\item{[Ba]} {\sc P.F. Baum}  {\tensl On the cohomology of homogeneous
spaces}, Topology 7 (1968), 15--38.

\item{[BB]} {\sc P.F. Baum, W. Browder} {\tensl The cohomology of
quotients of classical groups} Topology 3 (1965), 305--336.

\item{[Bu]} {\sc D. Burns}, {\tensl Harmonic maps from $\CP 1$ to $\CP n$},
in harmonic maps, Proc.,  New Orleans, 1980; LNM 949, Springer-Verlag,
Berlin/New York, 1982.

\item{[CCMM]}  {\sc F.R. Cohen, R.L. Cohen, B.M. Mann,
R.J. Milgram}, {\tensl The topology of rational functions and divisors
of surfaces}, Acta Math. 166(1991), 163--221.

\item{[CS]} {\sc R.L. Cohen, D. Shimamato}, {\tensl Rational functions,
labeled configurations, and Hilbert schemes}, J. London Math. Soc.
43(1991), 509--528.

\item{[C]} {\sc T.A. Crawford}, {\tensl Full holomorphic maps from the
Riemann sphere to complex projective space}, Jour. Diff. Geom. 38(1993),
161--189.

\item{[DZ]} {\sc A.M. Din, W.J. Zakrzewski}, {\tensl General
classical solutions in the $\CP{n-1}$ model}, Nuclear Phys. B 174
(1980), 397--406.

\item{[EW1]} {\sc J. Eells, J.C. Wood}, {\tensl Restrictions on
harmonic maps of surfaces},  Topology, 17(1976), 263--266.

\item{[EW2]} {\sc J. Eells, J.C. Wood}, {\tensl Harmonic maps from
surfaces to complex projective spaces}, Adv. in Math. 49(1983),
217--263.

\item{[FGKO]} {\sc M. Furuta, M.A. Guest, M. Kotani, Y. Ohnita}, {\tensl
On the fundamental group of the space of harmonic 2-spheres in the
$n$-sphere},  Math. Zeit. 215(1994), 5003--518.

\item{[GH]} {\sc P. Griffiths, J. Harris}, {\tensl Principles
of Algebraic Geometry},  Wiley, New York, 1978.

\item{[G]} {\sc M.A. Guest}, {\tensl Harmonic two-spheres in complex
projective space and some open problems}, Exp. Math., 10(1992) 61--87.

\item{[GO]} {\sc M.A. Guest, Y. Ohnita}, {\tensl Group actions and
deformations for harmonic maps}, J. Math. Soc. Japan, 45(1993)
671--704.

\item{[H]} {\sc J. Havlicek}, {\tensl The cohomology of holomorphic
self-maps of the Riemann sphere}, to appear Math. Zeit.

\item{[K]}  {\sc M. Kotani}, {\tensl Connectedness of the space of
minimal 2-spheres in $S^2m(1)$}, Proc. Amer. Math. Soc., 120(1994)
803--810.

\item{[La]} {\sc H.B. Lawson}, {\tensl La Classification des 2-sph\`eres
minimales dans l'espace projectif complex}, Ast\'erisque, 154--5(1987),
131--149.

\item{[LW]} {\sc L. Lemaire, J.C. Wood}, {\tensl On the space
of harmonic 2-spheres in $\CP 2$}, preprint (dg-ga/9510003).

\item{[Lo]} {\sc B. Loo}, {\tensl The space of harmonic maps of $S^2$
into $S^4$}, Trans. Amer. Math. Soc. 313(1989), 81--103.

\item{[S]} {\sc G. Segal}, {\tensl The topology of spaces of rational
functions}, Acta Math., 143(1979), 39--72.

\item{[V]} {\sc J.L. Verdier}, {\tensl Two dimensional $\grs$-models and
harmonic maps from $S^2$ to $S^{2n}$}, Lecture Notes in Physics 180,
Springer (Berlin), 1983, 136--141.

\item{[Woo]} {\sc G. Woo}, {\tensl Pseudo-particle configurations in
two-dimensional ferromagnets}, J. Math. Phys., 18(1977), 1264.

\bigskip
{\ninerm\halign{\indent#\hfill\cr
T.A. Crawford\cr
Fields Instutute for Research in Mathematical Sciences\cr
email: acrawfor@fields.utoronto.ca\cr}}

\vfill\eject
\bye